\documentclass[12pt]{amsart}
\usepackage{geometry} 
\usepackage{subfigure}
\usepackage{multicol}
\usepackage{graphicx}
\usepackage{amssymb,latexsym,amsmath,mathrsfs} 

\title[Normoglycemia and dysglycemia indices]{A validation study of normoglycemia and dysglycemia indices as a diabetes risk model}
\author{Paola Vargas}
\author{Miguel Angel Moreles}
\author{Joaquin Pe\~{n}a}
\address{Centro de Investigaci\'{o}n en Matem\'{a}ticas, Jalisco s/n, Valenciana,
Guanajuato, GTO 36240, Mexico, moreles@cimat.mx}

\author{ Adriana Monroy}
\address{Mexico General Hospital,
Mexico City, Mexico}
\date{} 

\begin{document}

\begin{abstract}
In this work, we test the performance of Peak glucose concentration ($A$) and average of glucose removal rates ($\alpha$), as normoglycemia and dysglycemia indices on a population monitored at the Mexico General Hospital between the years 2017 - 2019. A total of 1911 volunteer patients at the Mexico General Hospital are considered. 1282 female patients age ranging from 17 to 80 years old, and 629 male patients age ranging from 18 to 79 years old. For each volunteer, OGTT data is gathered and indices are estimated in Ackerman's model. A binary separation of normoglycemic and disglycemic patients using a Support Vector Machine with a linear kernel is carried out. Classification indices are successful for 83\%. Population clusters on diabetic conditions and progression from Normoglycemic to T2DM may be concluded. The classification indices, $A$ and $\alpha$ may be regarded as patient's indices and used to detect diabetes risk. Also, criteria for the applicability of glucose-insulin regulation models are introduced. The performance of Ackerman's model is shown.
\end{abstract}

\maketitle
\tableofcontents

\section{Introduction}

The diabetes situation is alarming worldwide, thousands of lives are lost and billions of dollars are spent on fighting the epidemic. In Mexico the trend of people having diabetes is of great concern. According to \cite{Levaillant_etal} 10.4\% of population had diabetes in 2016, compared with 7\% of the population in 2006. Comparative percentages are presented in Basto-Abreu et al \cite{Basto-Abreu_etal}, they report 16.8\%  and 15.7\% of the population with diabetes in 2018 and 2020 respectively. Making matters worse, a high percentage of adults are unaware of their risk to develop diabetes.
Characterizing such a risk is an important line of research, literature is too vast for an exhaustive review. We shall refer only to Gerstein et al. \cite{Gerstein_etal} and Tab\'{a}k et al \cite{Tabak_etal}. See also Hou et al \cite{Hou_etal}. 

In Gerstein et al. \cite{Gerstein_etal}, an estimate of the risk of progression to diabetes is presented. The estimates are based on a systematic literature overview of works that report data studies of individuals with impaired fasting glucose (IFG) or impaired glucose tolerance (IGT). These conditions as diagnosed by an oral glucose tolerance test (OGTT).  It is common for a subject to develop both conditions (IGT-IFT). A related study is that of Ardahanli \& Celik \cite{Ardahanli_Celik}, they argue that the progression to diabetes is expected for a significant percentage of people,  5-10\%. 

In Taba\'{a}k et al \cite{Tabak_etal}, an in-depth discussion on prediabetes is carried out. The discussion on risk prediction based on indices, associated with diabetes prediction models, aiding to identify individuals at risk of developing diabetes is of particular interest. A selected list of prediction models is presented therein. 

Prediabetes is typically defined as blood glucose concentrations higher than normal, but lower than diabetes
thresholds. It is well known that many people with prediabetes do not necessarily progress to diabetes, and the risk  
need not be different from people with other diabetes risk factors. Consequently, alternative points of view are welcome. In this work, instead of using the term prediabetes, we shall refer explicitly to dysglycemia conditions, IFG, IGT, and Type 2 Diabetes Mellitus (T2DM). A person with Normal Glucose Tolerance (NGT) shall be referred to as normoglycemic.  These categories are according to OGTT data using the American Diabetes Association (ADA) criteria, \cite{ADA}.

In line with Tab\'{a}k et al \cite{Tabak_etal}, a proof of concept study is presented in Vargas et al \cite{Vargas_etal}. Using OGTT data, patient indices are proposed to detect diabetes conditions, IGT, IFG,  Type 2 Diabetes Mellitus (T2DM), and Normal Glucose Tolerance (NGT). The proposed indices are peak glucose concentration ($A$), and the average of glucose removal rates ($\alpha$), two parameters in the glucose-insulin regulation model of Ackerman \cite{Ackermanetal1}. A binary classification learning technique is developed in terms of these indices. In a classical train-test split fashion, a population of 80 female subjects is used for training, and a mixed population, 24 males, and 33 females for testing. It results in 85\% of subjects being correctly classified. Also, it is observed a clockwise progression from normoglycemic to T2DM in the $A-\alpha$ plane. To be more conclusive, this observation requires to be tested on larger populations. The latter is the purpose of this work, we consider a population of 1911 subjects monitored at the Mexico General Hospital between the years 2017 - 2019.  Following the methodology introduced in Vargas et al \cite{Vargas_etal}, $A$ and  $\alpha$  are estimated for each patient in Ackerman's model for the corresponding OGTT data.  We shall argue that the clockwise progression from normoglycemic to T2DM remains valid, supporting the potential of these indices to detect diabetes risk.\\

It is apparent that conclusions depend upon the underlying glucose-insulin model.  We provide criteria for applicability and argue that, in spite of its simplicity, Ackerman's model is satisfactory. For better results, generalizations and /or modifications are to be sought.

The motivation for our research is to extract additional information from the data collected during an OGTT. This is by no means new. In \cite{Ismailetal}, the shape of the glucose concentration curve is used to predict the risk for type 1 diabetes. It is also associated with cardiovascular disease in  \cite{Hulmanetal}. Gaining insight into data is common practice. In the general context of endocrine diseases, Salazar et al \cite{Salazar_etal} make a strong case for data mining in classification.

Our work is not a classical statistical study of a control group. In the jargon of Data Analysis, our approach meets several of the six divisions of Data Science, Donoho \cite{Donoho}. In particular, data exploration and preparation, data modeling, data representation, and transformation.

\section{Materials and Methods}

\subsection{Oral Glucose Tolerance Test data and ADA classification}
A total of 1911 volunteer patients participated in a monitoring study at the Mexico General Hospital between the years 2017 - 2019. 
Patients came from several states in Mexico. 1282 female patients age ranging from 17 to 80 years old, and 629 male patients age ranging from 18 to 79 years old. 

The Ethical Committee of the Mexico General Hospital approved the study protocol and informed written consent was obtained from each volunteer.

Volunteers underwent a 2h OGTT to determine normoglycemia (Normal Glucose Tolerance, NGT) or dysglycemia. The latter, any of the diabetic conditions:  Impaired Fasting Glucose (IFG), Impaired Glucose Tolerance (IGT), and Type 2 Diabetes Mellitus (T2DM).  The American Diabetes Association (ADA) criteria were used for classification, \cite{ADA}. Namely,

\begin{itemize}
\item Diabetes Mellitus: fasting plasma glucose $\geq$ 7 mmol/l (126 mg/dl) or 2 h post-load plasma glucose $\geq$ 11.1 mmol/l (200 mg/dl). 

\item IFG: fasting plasma glucose $\geq$ 5.6 to 6.9 mmol/l (100-125 mg/dl),

\item  IGT: as 2 h plasma glucose $\geq$ 7.8 to 11 mmol/l (140-199 mg/dl). 
\end{itemize}

As a result, the criteria yield 1084 NGT, 147 IFG, 319 IGT, 156 IFG-IGT, and 205 T2DM. 

\subsection{Indices for normoglycemia and dysglycemia}

For a given subject, let $G(t)$ be the plasma glucose concentration in the blood at time $t$.  
The graph of the function $G(t)$ is the so called, OGTT curve (pattern). \\

It is assumed that after fasting, the patient's concentration has stabilized to $G_0$. 
 
After the glucose load has been absorbed, Ackerman's model  \cite{Ackermanetal1} yields
 \[
G(t) = G_0 +  Ae^{-\alpha t} \cos(\omega t - \delta).
 \]
 
Consequently, if the parameters $A, \alpha,\omega,\delta$ are known, the corresponding OGTT pattern is obtained.
In Vargas et al \cite{Vargas_etal}, these parameters are estimated using a Bayesian approach.  In practice, a sample of the posterior density is constructed and point estimators are computed. In this case, the Maximum A Posteriori (MAP) estimator.

The biophysically meaningful parameters proposed  for classification are:

\begin{itemize} 
\item $A$, \emph{peak glucose concentration},

\item $\alpha$, \emph{average of glucose removal rates},
\end{itemize}

A classical training test strategy was carried out. The training set is comprised of 80 (female) subjects. These are divided in two classes, normoglycemic (+1) and disglycemic  (-1). For patient $j$, the parameters $A^j$ and $\alpha^j$ are considered for binary classification in the $A$-$\alpha$ plane.  A linear separator is found by the soft margin support vector machine methodology, see \cite{JWHT}.\\

The test set is an independent mixed population of comparable size:

\begin{itemize}

\item 24 Males: 15 NGT, 3 IGT, 2 IFG-IGT, 4 T2DM,

\item 33 Females: 17 NGT, 1 IFG, 11 IGT, 2 IFG-IGT, 2 T2DM.

\end{itemize}

The resulting classification was successful for $91\%$ of males and $87\%$ of females. A relevant conclusion is the apparent \emph{clockwise} progression from normoglycemic to Type 2 Diabetes Mellitus. This fact is illustrated  in the 
$A-\alpha$ plane with five patients, one for each condition, see Figure 1. It is argued that there is a progression from normoglycemic to disglycemic at the separating line. 

\begin{figure}[htbp]
\begin{center}
  \includegraphics[width=3.8in]{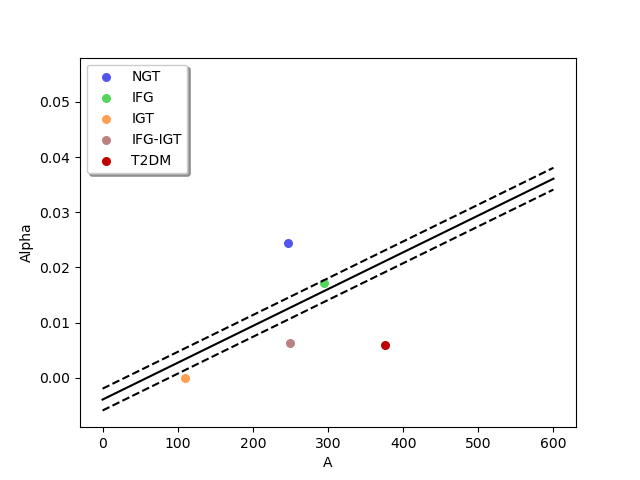}
\caption{Clockwise progression from NGT to T2DM.}
\label{default}
\end{center}
\end{figure}

It is concluded that peak glucose concentration and average of glucose removal rates are potentially normoglycemia and dysglycemia indices. 

\subsection{Criteria for applicability of Ackerman's model}

Let us denote the glucose concentration in minute $i$  by $G_i$. 

Given OGTT data for a subject, $G_0, G_{30},G_{60}, G_{90},G_{120},$ the parameters $A, \alpha,\omega,\delta$ are estimated. 

By means of Ackerman's model, predicted concentrations are obtained. Namely, $G_{i}^{pred}$ at minute $i$.

It is customary to consider a model to be appropriate if there is a good fit for the data.  We shall use the latter as part of the criteria. 

Consider the total absolute error between the patient's glucose concentration data and the predicted concentration by the estimated parameters in Ackerman's model.  Namely,
 $$\text{Error}_{abs} = \frac{1}{5}\sum_{i\in \{0,30,60,90,120\}}|G_{i} - G_{i}^{pred}|.$$
 
In the context of the OGTT, there are some biomedical requirements. The OGTT curve is wave-like, but oscillations in the 120 minutes time interval are not arbitrary. In Kanauchi et al  \cite{Kanauchi_etal} it is argued that for an OGTT curve to be admissible, at most two oscillations are realistic and maximum concentration is attained. Also, a quantitative criterion is indicated in Tschritter et al \cite{Tschritter_etal}. It is proposed that the difference between glucose concentrations of any given patient at $90$ and $120$ minutes, $G_{90}$ and $G_{120}$ should exceed 4.5 mg/dl (0.25 mmol/l).

Combining these biomedical conditions with data fitness, we are led to the following criteria.

Given OGTT data for a single subject, we consider that Ackerman's model applies to this data if  $\omega < 0.09$, and one of the following conditions is satisfied:

\begin{enumerate}
    \item Data is fitted accurately, $\text{Error}_{abs} < \text{4.5}.$
     
    \item The shape of curve is admissible, $\vert G_{90} - G_{120}\vert < 4.5$ mg/dl and $\text{Error}_{abs} < 5$. 
    
    \item The shape of the curve is admissible, $\vert G_{90} - G_{120}\vert > 4.5$ mg/dl; $\text{Error}_{abs} < \text{7.5}.$
\end{enumerate}

\medskip

With $\omega < 0.09$ oscillations of the OGTT curve are limited to two. In Condition (3), the fit is not stringent, but both biomedical requirements are valid.
 
\section{Results}

\subsection{Total population}
Parameter estimation is carried out for the 1911 OGTT data cases, and the constructed linear separator is now tested on this heterogeneous population.

\begin{figure}[!ht]
    \begin{center}
        \includegraphics[scale=0.6]{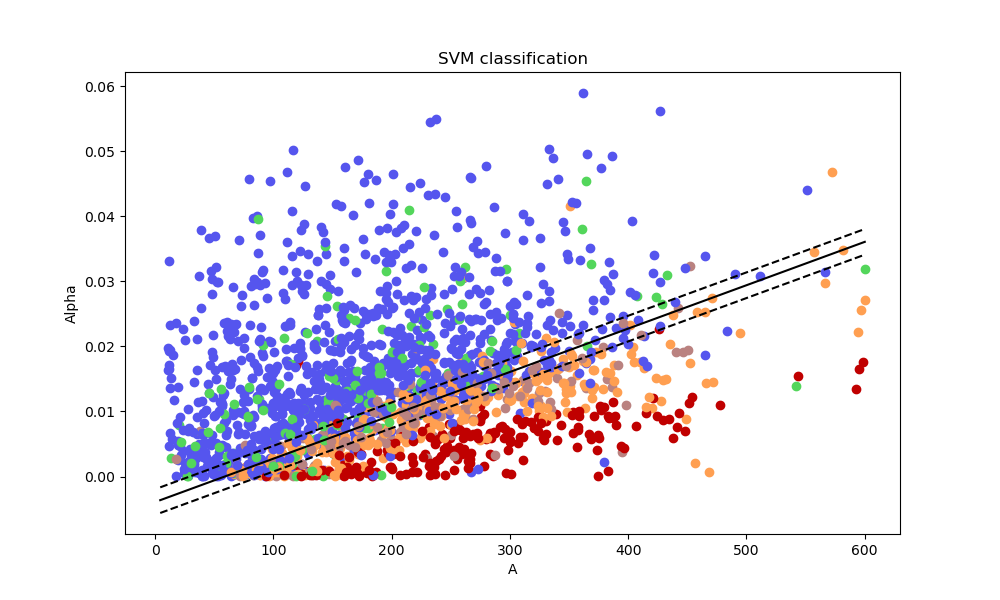}
        \caption{Classification diagram for the 1911 OGTT data cases.}
        \label{f_ttt_pc}
    \end{center}
\end{figure}

The population clusters are as follows: 1084 NGT, 147 IFG, 319 IGT, 156 IFG-IGT, and 205 T2DM. 
The proposal is to consider $A$ and $\alpha$ as classification indices. The estimated parameters for the total population are shown in Figure \ref{f_ttt_pc}. Classification is 82.26 \% successful.

As in the proof of concept, there are three clusters preserving the clockwise progression:
\begin{itemize}
\item Normoglycemic and IFG patients,

\item IGT patients, or patients with both conditions IFG and IGT

\item Type 2 diabetes mellitus patients
\end{itemize}

Classification of a T2DM patient as normoglycemic is the most costly error.  As shown in Figure \ref{f_ttt_pc}, most of the subjects with T2DM condition are below the separating line. 

In terms of each subgroup. Correct classification is as follows:\\

\begin{itemize}
 \item NGT, 90.12\%
 \item IFG,  21.08\%
 \item IGT, 76.8\%
 \item IFG-IGT, 75\%
 \item DMT2, 98.53\%
 \end{itemize}

The IFG patients cluster with normoglycemic may be due to the fact that de glucose pattern of both populations is similar at the end of the OGTT. Also, IFG is the mildest condition.

\subsection{Population filtering}

Parameter estimation is carried out and the applicability of Ackerman's model is tested. It is found that the model is valid for 1302 subjects, $68.12\%$ of the population.

In Figure 3 we show some examples of OGTT curves where Ackerman's model is valid.

\begin{figure}[!ht]
    \begin{center}
        $\begin{array}{ccc}
             \includegraphics[width=1.8in]{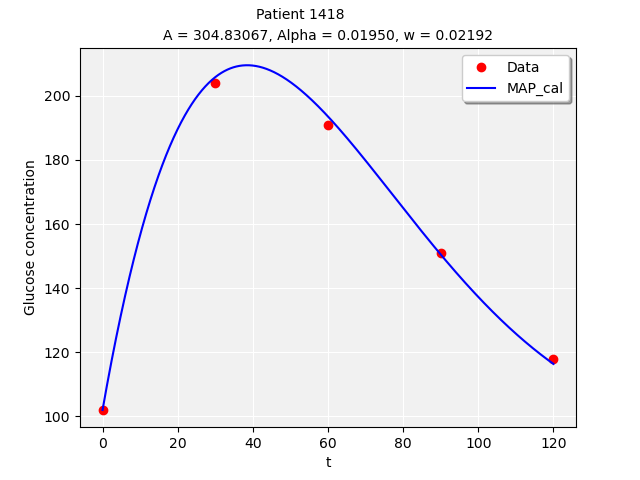} & \includegraphics[width=1.8in]{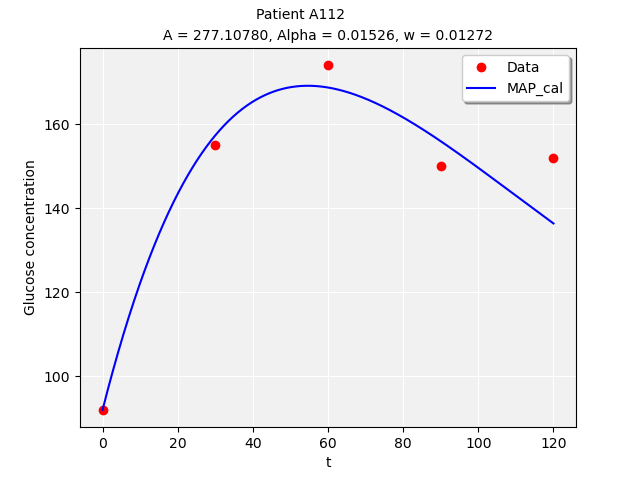}  & \includegraphics[width=1.8in]{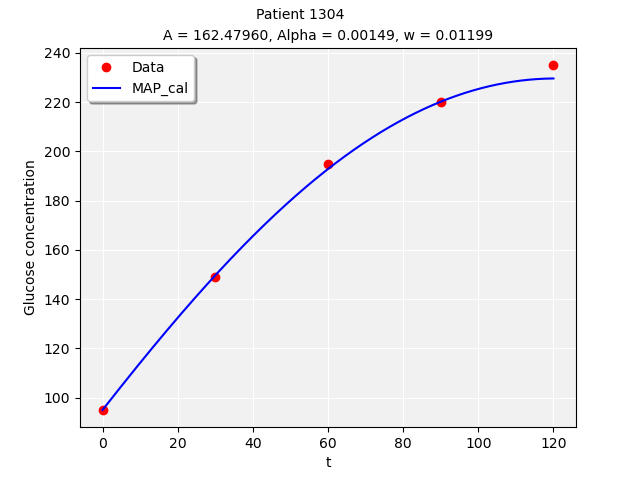} 
             \\
            \text{(a) Condition (1) is valid} & \text{(b) Condition (2) is valid} & \text{(c) Condition (3) is valid}
        \end{array}$
    \end{center}
    \caption{OGTT curves represented by Ackerman's model, $\omega<0.09$.}
    \label{H_C}
\end{figure}

In Figure 4 some examples where the model fails. The last figure illustrates that fitness need not be enough.

\begin{figure}[!ht]
    \begin{center}
        $\begin{array}{ccc}
           \includegraphics[width=1.8in]{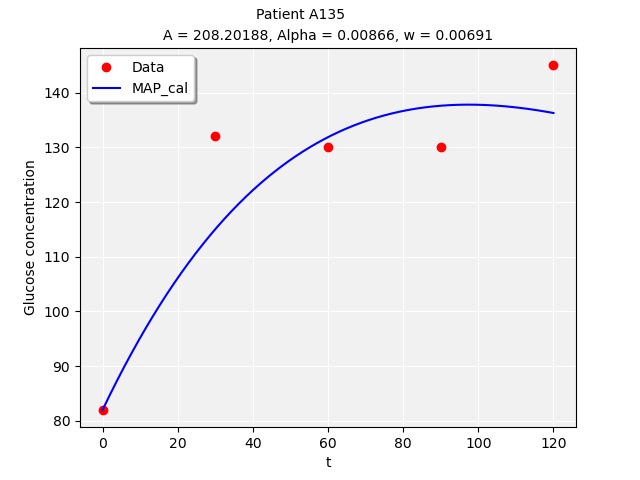} &             \includegraphics[width=1.8in]{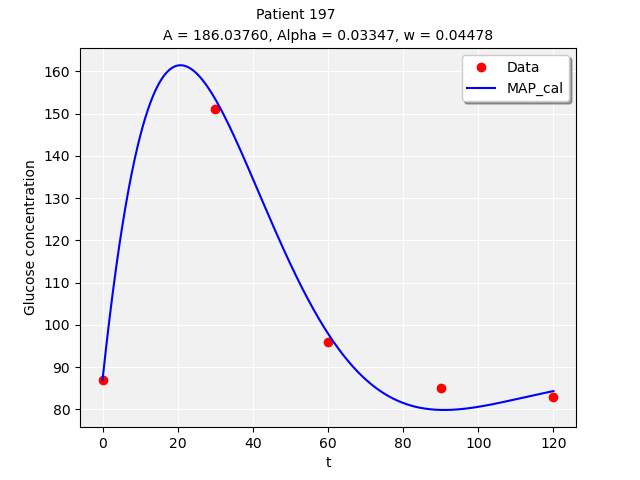} &
            \includegraphics[width=1.8in]{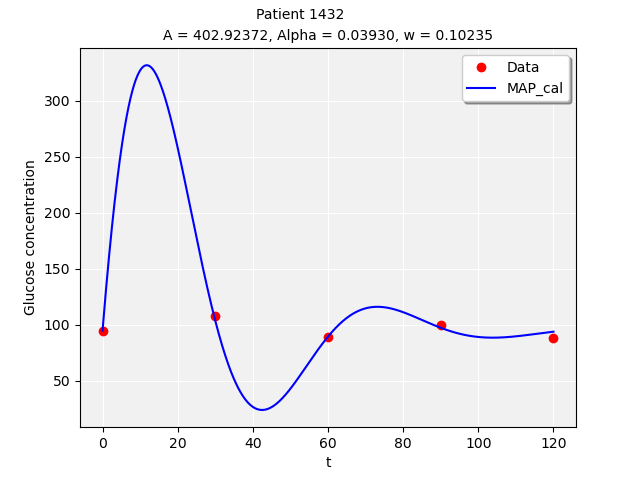}\\
            \text{(a) Condition (1) fails} & \text{(b) Condition (2) fails } & \text{(c) $\omega>0.09$ }
        \end{array}$
    \end{center}
    \caption{OGTT curves not represented by Ackerman's model}
    \label{H_C}
\end{figure}

The resulting group clusters are as follows: 687 NGT, 102 IFG, 186 IGT, 106 IFG-IGT, and 129 T2DM.
The estimated parameters for the filtered population are shown in Figure \ref{default} in the $A-\alpha$ plane. Classification is successful for 83.56\%
Results are consistent with the full population.

\begin{figure}[htbp]
\begin{center}
  \includegraphics[width=3.8in]{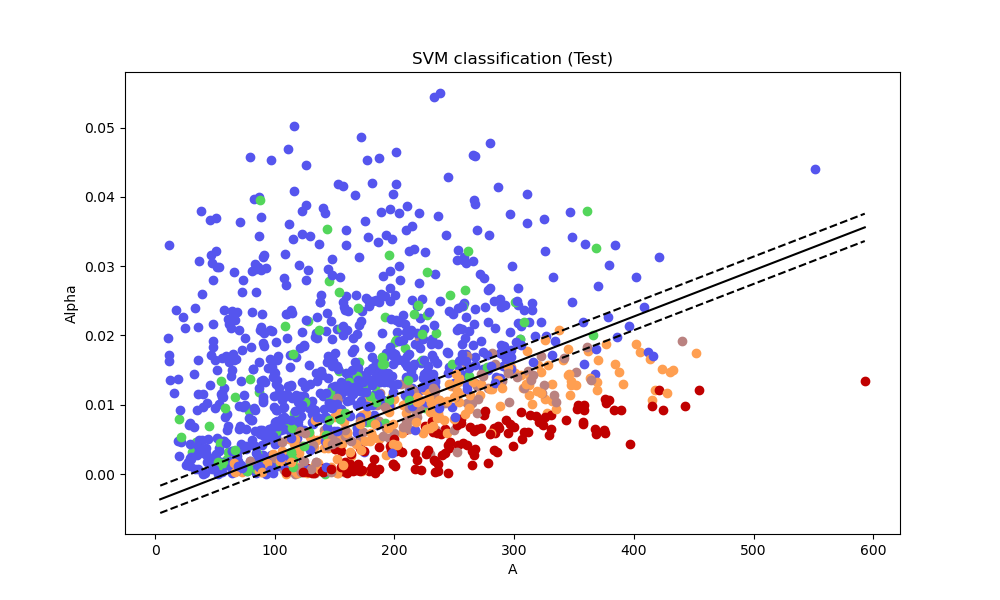}
\caption{Classification for the filtered population.}
\label{default}
\end{center}
\end{figure}

\section{Discussion and conclusions}

The main objective of this work was to test the methodology developed in Vargas et al \cite{Vargas_etal}, on the 
heterogeneous population monitored in the Mexico General Hospital. We contend that the results obtained on the heterogeneous population, make a case for the conclusions in Vargas et al \cite{Vargas_etal}:

\begin{itemize}
\item The parameters $A$, $\alpha$ can be regarded as patient indices to determine normoglycemia and dysglycemia

\item There is a linear strip separating normoglycemic from disglycemic patients in the $A-\alpha$ plane

\item There is a clockwise progression from normoglycemic to dysglicemic in the $A-\alpha$ plane

\item A patient with point indices in a neighborhood of the linear strip may be regarded as a sign of diabetic risk
\end{itemize}

From the modeling perspective, we proposed criteria for the applicability of Ackerman's model to describe the glucose-insulin regulatory system during an OGTT. Notice that the criteria do not depend on the parameter estimation method.\\

The Geometric criterion proposed in Kanauchi \cite{Kanauchi_etal} was realized by limiting the value of $\omega$. The criterion can be enforced by the analog using the period of the curve, 
\[
T=\frac{2\pi}{\omega}.
\] 
Hence $\omega<0.09$ can be replaced by 
\[
T>\frac{2\pi}{0.09}\approx  70 \text{ minutes.}
\]

Consequently, it is apparent that the resulting new criteria can be applied to any OGTT curve model and any parameter estimation methodology. 

It is remarkable that in all its simplicity, Ackerman's model describes correctly almost 70\% of the glucose patterns. We stress again that the population under study was not a control group. The problem was approached as data exploration, data representation and modeling. The obtained results suggest seeking generalizations and/or modifications of Ackerman's model for glucose-insulin regulation. 

In the database, there are a few patients with two or three Oral Glucose Tolerance Tests.  In a supervised patient case, the $A$-$\alpha$ indices are potentially useful to monitor diabetic conditions. In a preliminary fashion,  let us show results for a patient with three sets of data ordered chronologically. Figure 5 shows the evolution of the $A,\alpha$ indices.

 \begin{figure}[!ht] 
      \begin{center}
      \includegraphics[scale=0.48]{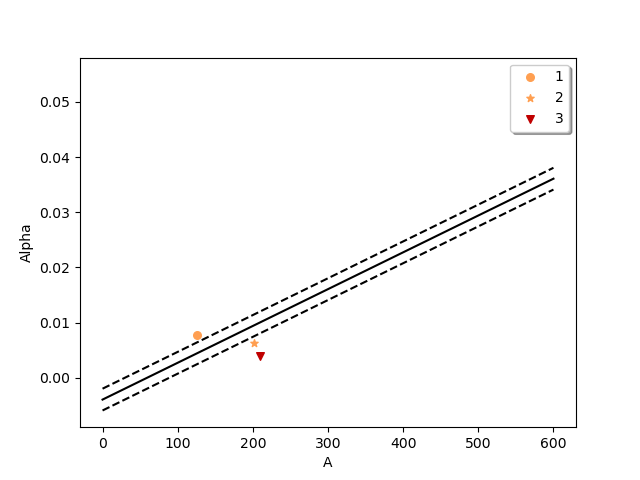} 
       \end{center}
            \caption{\scriptsize{Evolution of diabetic condition.}}
 \end{figure}

It is apparent the progression from IGT to IGT-IFG to T2DM. A thorough monitoring study is required to test this hypothesis.

An exploration of all these issues is left for future works.

\end{document}